\documentstyle[12pt]{article}

\def\Empty{}
\catcode`\@=11
\def\section{\@startsection {section}{1}{\z@}{-3.5ex plus-1ex minus
    -.2ex}{2.3ex plus.2ex}{\large\bf}}
\def\subsection{\@startsection{subsection}{2}{\z@}{-3.25ex plus-1ex
    minus-.2ex}{1.5ex plus.2ex}{\large\bf}}
\def\subsubsection{\@startsection{subsubsection}{3}{\z@}{-3.25ex plus
 -1ex minus-.2ex}{1.5ex plus.2ex}{\normalsize\bf}}

\def\eqalign#1{\,\vcenter{\openup\jot\m@th
  \ialign{\strut\hfil$\displaystyle{##}$&$\displaystyle{{}##}$\hfil
        \crcr#1\crcr}}\,}

\def\mydesc{\list{}{\labelwidth\z@ \itemindent-\leftmargin
\listparindent 1.5em
\let\makelabel\descriptionlabel}}

\def\fnum@figure{{\small Figure \thefigure}}
\def\fakefigure{\def\@captype{figure}}

\long\def\@makecaption#1#2{
    \vskip 10pt
    \def\FCap{#2} \def\NoCap{\ignorespaces}
    \ifx \FCap\NoCap
       \setbox\@tempboxa\hbox{#1}  % This is to avoid the damn colon.
      \else
       \setbox\@tempboxa\hbox{#1: \small \it #2}
    \fi
    \ifdim \wd\@tempboxa >\hsize   % IF longer than one line:
        \unhbox\@tempboxa\par      %   THEN set as ordinary paragraph.
      \else                        %   ELSE  center.
        \hbox to\hsize{\hfil\box\@tempboxa\hfil}
    \fi}

\def\@oddhead{\hbox{}\rightmark \hfil \rm\thepage}% Right heading.
\def\sectionmark#1{\markright {\sc{\ifnum \c@secnumdepth >\z@
      \S\thesection.\hskip 1em\relax \fi #1}}}

\catcode`\@=12

\def\oplabel#1{
  \def\OpArg{#1} \ifx \OpArg\Empty {} \else
        \label{#1}
  \fi}

\def\MakeStEnv#1{
  \newenvironment{#1}[2]{
  \begin{#1St} \oplabel{##1}%
  \global\def\CrntSt{\thetheoremSt}%
  {\rm ##2}%
}{
  \end{#1St} }
}
\MakeStEnv{theorem}
\MakeStEnv{corollary}
\MakeStEnv{proposition}
\MakeStEnv{lemma}
\MakeStEnv{conjecture}

\newenvironment{proof}[1]{
  \def\PfArg{#1}
  \ifx\PfArg\Empty
        \edef\PfArg{\CrntSt}  \fi
 \startproof{\PfArg}%
}{
  \finishproof{\PfArg}
}
\newcommand{\startproof}[1]{
  \medbreak\mbox{}
  {\it Proof of #1:}%
}
\newcommand{\finishproof}[1]{
  \def\FPArg{#1}
  \ifx\FPArg\Empty
        \def\FPArg{\CrntSt}  \fi
  \smallbreak\noindent\makebox[\textwidth]{\hfill\fbox{\FPArg}}
  \medbreak\noindent
}

\newcommand{\marginwrite}[1]{}
\def\abs#1{\vert#1\vert}     % absolute value
\def\del{\partial}           % boundary symbol
\def\eq{\!=\!}               % improved spacing for equals sign
\def\frac#1#2{{{\displaystyle #1}\over{\displaystyle #2}}}   % fractions
\def\rel{\hbox{\rm \ rel\ }} % rel in math mode
\def\set#1{\{ #1 \}}         % set braces  {  }
\def\text#1{\hbox{{\rm#1}}}  % creates an hbox in roman font
\def\Diff{\hbox{\it Diff\/}}
\def\diff{\hbox{\it diff\/}}

\def\Imb{\hbox{\it Imb\/}}
\def\imb{\hbox{\it imb\/}}
\def\Isom{\hbox{\it Isom\/}}
\def\isom{\hbox{\it isom\/}}

\def\Norm{\hbox{\rm Norm}}
\def\norm{\hbox{\rm norm}}
\def\O{\text{O}}
\def\Out{\hbox{\rm Out}}
\def\SO{\text{SO}}
\def\Spin{\hbox{\rm Spin}}
\def\SU{\text{SU}}

\def\mapdown#1{\big\downarrow % makes a downarrow with a mathematics
            \rlap            % mode symbol to the right
            {\smash{$\vcenter       % (for use in the "diagram" template)
            {\hbox{$         %
            \scriptstyle#1   %
            $}}$}}}           %

\def\mapright#1{\smash{      % makes a long rightarrow with a mathematics
            \mathop          % mode symbol above it
            {\longrightarrow % (for use in the "diagram" template)
            }\limits^{#1}}}  %

\def\P{\hbox{\rm I\kern-1.8pt P}}            % projective plane numbers symbol

\def\R{\hbox{\rm I\kern-1.8pt R}}            % real numbers symbol

\def\Z{\hbox{\rm Z\kern-3.5pt Z}}            % integers symbol

\font\small=cmr10 scaled 833

\def\G{\Gamma}

\def\t{\tau}

\begin{document}

\title{The Generalized Smale Conjecture for 3-manifolds with genus 2
one-sided Heegaard splittings}
\author{Darryl McCullough and J. H. Rubinstein}
\date{{\footnotesize Department of Mathematics, University of Oklahoma,
Norman, OK 73019, USA}
\\
{\footnotesize Department of Mathematics, University of
Melbourne, Parkville, Victoria 3052, Australia}
\\
\bigskip
{\footnotesize Version of November 9, 1996}
}
\maketitle

\section[Introduction]{Introduction}
\label{intro}

The Smale Conjecture, proved by Hatcher \cite{H1}, asserts that if $M$
is the 3-sphere with the standard constant curvature metric, the
inclusion $\Isom(M)\rightarrow \Diff(M)$ from the isometry group to
the diffeomorphism group is a homotopy equivalence. The {\it
Generalized Smale Conjecture} asserts this whenever $M$ is a closed
3-manifold with a metric of constant positive curvature. All known
examples of closed 3-manifolds with finite fundamental group admit
such metrics, whose isometry groups are easily determined. For many
manifolds of constant positive curvature, $\Isom(M)\rightarrow
\Diff(M)$ is known to be bijective on the set of path components
\cite{B}, \cite{R-B}, {\it et.~al.}

The Generalized Smale Conjecture is an instance of the general
principle, first realized by Thurston, that 3-manifold topology is
profoundly affected by the existence and behavior of geometric
structures. In the positive curvature case, the Generalized Smale
Conjecture suggests that not only the topological structure but also
the group of all smooth automorphisms is controlled by the
geometry. For (compact) 3-manifolds whose interiors have constant
negative curvature and finite volume, the analogous expectation holds
true, at least when the manifolds are sufficiently large, since by
Mostow Rigidity, Waldhausen's Theorem, and work of Hatcher discussed
below, the composition $\Isom(M)\to \Out(\pi_1(M))\to \Diff(M)$ is a
homotopy equivalence. In contrast, when the manifold has interior of
constant negative curvature and infinite volume, or has constant zero
curvature, a diffeomorphism will not in general be isotopic to an
isometry (said differently, the diffeomorphism group may have more
components than the isometry group). Even in these cases, however,
Waldhausen's Theorem and Hatcher's work show that $\Isom(M)\to
\Diff(M)$ is always a homotopy equivalence when one restricts to the
connected components of the identity diffeomorphism.

Our approach to the Generalized Smale Conjecture, for the cases we
will consider here, is based on work of Hatcher as extended by Ivanov.
For sufficiently large $\P^2$-irreducible 3-manifolds, Hatcher
(\cite{H2}, combined with~\cite{H1}), extending earlier work of
Laudenbach \cite{L}, proved that the components of $\Diff(M \rel
\partial M)$ are contractible. The main part of the argument is to
show that the space of imbeddings of a 2-sided incompressible surface
$F$ that are disjoint from a parallel copy of $F$ is a deformation
retract of the space of all imbeddings of~$F$ (isotopic to the
inclusion relative to $\partial F$).

For manifolds that are not sufficiently large and therefore do not
contain a 2-sided incompressible surface, one may try to use a 1-sided
incompressible surface instead. If $M$ is orientable, irreducible, and
not sufficiently large, and contains a 1-sided incompressible surface
$K$, then by theorem~4 of \cite{R1}, $M-K$ is an open handlebody. When
the complement of a 1-sided surface $K$ in $M$ is an open handlebody,
we say that $(M,K)$ is a {\it 1-sided Heegaard splitting}. The {\it
genus} of the splitting is the (nonorientable) genus of $K$. Ivanov
\cite{I1}, \cite{I2} extended Hatcher's results to a restricted class of
3-manifolds with genus 2 1-sided Heegaard splittings, enabling  him to
determine the homotopy type of their diffeomorphism groups. In
particular, the Generalized Smale Conjecture is implied for these
manifolds.

The manifolds considered by Ivanov are those for which no Seifert
fibering is nonsingular on the complement of any vertical Klein
bottle. The remaining cases are:

\begin{enumerate}
\item[I.] Those for which every Seifert fibering nonsingular on the
complement of $K$ restricts to the ``meridinal'' (nonsingular)
fibering of $K$. These have binary dihedral fundamental groups.

\item[II.] Those for which every Seifert fibering nonsingular on the
complement of $K$ restricts to the ``longitudinal'' (two exceptional
fibers) fibering of $K$. These are the lens spaces
$L(4n,2n-1),\;n\geq2$.

\item[III.] The lens space $L(4,1)$ (which admits both meridinal and
longitudinal fiberings nonsingular on the complement of $K$).
\end{enumerate}

\noindent In this paper, we prove the Conjecture for cases I and II.
More precise statements and an outline of the argument will be given
in section~\ref{results}, after fixing notation.

Ivanov's announced results were used in \cite{F-W} to construct
examples of homeomorphisms of reducible 3-manifolds that are homotopic
but not isotopic. Our results show that their construction applies to
a larger class of 3-manifolds. In \cite{P-R}, our work was applied to
the classification problem for 3-manifolds which have metrics of
positive Ricci curvature and universal cover $S^3$.

The Generalized Smale Conjecture has attracted the interest of
physicists studying the theory of quantum gravity. Certain physical
configuration spaces can be realized as the quotient space of a
principal $\Diff_1(M,x_0)$-bundle with contractible total space, where
$\Diff_1(M,x_0)$ denotes the subspace of $\Diff(M,x_0)$ that induce
the identity on the tangent space to $M$ at $x_0$. (This group is
homotopy equivalent to $\Diff(M\# D^3\rel \del D^3)$.) Consequently
the loop space of the configuration space is weakly homotopy
equivalent to $\Diff_1(M,x_0)$. Physical significance of
$\pi_0(\Diff(M))$ for quantum gravity was first pointed out in
\cite{F-S}. See also \cite{ABBJRS}, \cite{G-L}, \cite{I}, \cite{S3},
\cite{W2}. The significance of some higher homotopy groups of
$\Diff(M)$ is examined in \cite{Giulini}.

An earlier version of this paper was circulated in the late 1980's,
and indeed has been cited a number of times in the scientific
literature. In this new version, the essential mathematical content is
unchanged, but a considerable amount of detail has been added. Also,
various ``folk'' theorems about fibrations of spaces of
diffeomorphisms and imbeddings, heavily used in our arguments, have
been put on firm ground by the work in~\cite{Kalliongis-McCullough}.

\section[Notation and statement of results]
{Notation and statement of results}
\label{results}

Let $K$ be a Klein bottle and let $P$ be the orientable $I$-bundle
over $K$ with boundary the torus $T$.  Let $R$ be a solid torus
containing a meridinal 2-disc whose boundary $C$ is a meridinal circle
lying in $\del R$. Fix a presentation $\pi_1(K)\eq
\langle\,a,\,b\,\,\vert\,\,b^{-1}ab\eq a^{-1}\, \rangle$. The four
homotopy classes of essential simple closed curves on $K$ are $b$,
$ab$, $a$, and $b^2$, with $b$ and $ab$ orientation-reversing and $a$
and $b^2$ orientation-preserving. The free abelian group
$\pi_1(\partial P)$ is generated by (loops homotopic in $P$ to) $a$
and $b^2$. For a pair $(m,n)$ of relatively prime integers, the
3-manifold $M(m,n)$ is formed by identifying $\partial R$ and
$\partial P$ in such a way that $C$ is attached along a simple closed
curve representing the element $a^mb^{2n}$. Since $M(-m,n)\eq M(m,n)$ and
$M(-m,-n)\eq M(m,n)$, we can and always will assume that both $m$ and $n$
are positive. The fundamental group of $M(m,n)$ has presentation
$\langle\,a,\,b\,\,\vert\,\,b^{-1}ab\eq a^{-1},\, a^mb^{2n}\eq
1\,\rangle$.

From \cite{R2} we have the following facts. If $m\eq 1$, then $M(1,n)$ is
the lens space $L(4n,2n-1)$. Suppose that $m\neq1$. If $m$ is even,
then $\pi_1(M(m,n))\cong D^*_{4m}\times C_n$, the direct product of
the binary dihedral group $D^*_{4m}\eq\langle
x,y\,\vert\,x^2\eq (xy)^2\eq y^{2m}\,\rangle$ and a cyclic group. Finally,
if $m>1$ is odd, write $4n\eq 2^kn_1$ where $n_1$ is odd. Then
$\pi_1(M(m,n))\cong D_{(2^k,m)}\times C_{n_1}$ where $D_{(2^k,m)}$ is
the generalized dihedral group
$\langle\,x,y\,\vert\,x^{2^k}\eq y^m\eq 1,\,x^{-1}yx\eq y^{-1}\,\rangle$. Note
that when $m$ is odd, there is an isomorphism from $D^*_{4m}$ to
$D_{(4,m)}$ given by sending $x$ to $x$ and $y$ to~$yx^2$.

In \cite{I1}, it was announced that for $n\neq1$, the inclusion
$\hbox{\it Isom}(M(m,n))\rightarrow\Diff(M(m,n))$ is a homotopy
equivalence, and a detailed proof for the case $m\neq1$ and $n\neq1$
was given in \cite{I2}. In the remaining sections of this paper, we
will prove this result for the cases where $m\eq1$ or $n\eq1$, except
for that of $m\eq n\eq 1$, where $M(1,1)\eq L(4,1)$. In
section~\ref{isometry}, we compute the isometry groups of the
$M(m,n)$. In section~\ref{mainthms}, we reduce the Conjecture in the
case of the $M(m,n)$ to proving that the inclusion from the space of
fiber-preserving imbeddings of $K$ into $M$ to the space of all
imbeddings is a weak homotopy equivalence on the connected components
of the inclusion. That is, a parameterized family $\set{K_t}_{t\in
D^k}$ of imbedded Klein bottles isotopic to $K$ can be deformed to a
fiber-preserving family. The proof of this assertion occupies the
final two sections. In contrast to the 2-sided case, one cannot
avoid parameters for which $K_t$ intersects $K$ non-transversely and
these configurations must be analyzed, but it is enough to consider
``generic position'' non-transverse configurations, as described in
\cite{I2}. For these we show in theorem~\ref{circles} that for each
parameter value $t$ one can find a concentric fibered torus $T_u$ in a
neighborhood of $K$ which meets $K_t$ transversely in circles that are
either inessential in $T_u$ or cover imbedded circles in $K_t$. In
section~\ref{parameterization}, we complete the argument, using the
methods of Hatcher to eliminate inessential intersections with the
concentric solid tori in $M-K$, then deforming the $K_t$ to be
fiber-preserving inductively over the skeleta of a triangulation of
the parameter space.

This program fails in a fundamental way in the case of $L(4,1)$
because no Seifert fibering of $L(4,1)$ is preserved by all
isometries. Below we will point out more precisely the steps where the
arguments break down.

\section[Calculation of isometry groups]
{Calculation of isometry groups}
\label{isometry}

The finite subgroups of $\SO(4)$ that act freely on $S^3$ were worked
out by Hopf and Seifert-Threlfall, and reformulated using quaternions
by Hattori. We review this as described in \cite{Wolf} (see especially
pp.~226-227).  There is a quotient map $F\colon\Spin(3)\times
\Spin(3)\to \SO(4)$, where $\Spin(3)\cong \SU(2)$ can be identified
with the group of unit quaternions, and $F(a,b)$ acts on $\R^4$ by
$F(a,b)(q)\eq aqb^{-1}$.  The kernel of $F$ is
$W\eq\set{(1,1),(-1,-1)}$.  The center of $\SO(4)$ has order~2, and is
generated by $[(1,-1)]$ which acts as the antipodal map on $S^3$.  Let
$G$ be a finite subgroup of $\SO(4)$ acting freely on $S^3$, and let
$M\eq S^3/G$. If $G$ has even order, then it must contain the
antipodal map.  As explained in $\cite{Orlik}$, the images of the two
$\Spin(3)$ factors in $\SO(4)$ can be described as the groups of
``right rotations'' and ``left rotations'', which commute and
intersect only in the antipodal map.  Let $G^*\eq F^{-1}(G)$, and let
$G_1$ and $G_2$ be the projections of $G^*$ into the factors of
$\Spin(3)\times\Spin(3)$. For $q\neq 0$, $F(a,b)(q)\eq q$ if and only
if $a\eq qbq^{-1}$, i.~e.~if and only if $a$ is conjugate to $b$ in
$\Spin(3)$.  Thus $G$ acts freely on $S^3$ if and only if $G^*$ has no
pair $(a,b)\notin W$ such that $a$ is conjugate to $b$ in $\Spin(3)$.
Upon detailed examination, this implies that at least one of $G_i$,
say $G_2$, is cyclic and hence is contained in a circle subgroup $S$
of $\Spin(3)$. Thus $F(S)$ is contained in the normalizer of $G$. This
implies that there is an action of $S^1$ by isometries on $M$, which
determines a Seifert fibering of~$M$. In~\cite{Orlik}, the explicit
imbeddings of the various $G$ into $\SO(4)$ are given, and we will
refer to these when we work out the isometry groups of some of the
quotient manifolds.

Since $\O(4)$ is the full group of isometries of $S^3$, the isometry
group of $M$ is the quotient $\Norm(G)/G$ where $\Norm(G)$ is the
normalizer of $G$ in $\O(4)$.  We are especially interested in
$\isom(M)$, the connected component of the identity in $\Isom(M)$.
Let $\norm(G^*)$ denote the connected component of the identity in the
normalizer of $G^*$ in $\Spin(3)\times \Spin(3)$.  Clearly the
connected component of the identity in $\Isom(M)$ is
$\norm(G^*)/(G^*\cap\norm(G^*))$.  For computing this, we observe that
$\norm(G^*)\eq\norm(G_1)\times\norm(G_2)$, where $\norm(G_i)$ denotes
the connected component of the identity in the normalizer of $G_i$ in
$\Spin(3)$.  When $G_i$ is cyclic of order~2, $\norm(G_i)\eq\Spin(3)$.
When $G_i$ is cyclic of order greater than~2, $\norm(G_i)$ is the
unique~$S^1$-subgroup of $\Spin(3)$ that contains~$G_i$.  When $G_i$
is noncyclic, $\norm(G_i)$ is just the identity element.

We now specialize to the manifolds $M(m,n)$ described in
section~\ref{results}. Let $G\eq \pi_1(M(m,n))$.

\smallskip

\noindent {\sl Case I.} $G$ is cyclic of order~4.

For $(a,b)\in G^*$, not both $a$ and $b$ can be of order~4, since they
would then be conjugate in $\Spin(3)$. Therefore one of the $G_i$ is
$C_4$ and the other is $C_2$, so $\norm(G^*)$ is $S^1\times
\Spin(3)$ and $\isom(L(4,1))$ is $S^1\times\SO(3)$.

\smallskip

\noindent {\sl Case II.} $G$ is cyclic of order~$4n$, $n\geq 2$.

Since the quotient is $L(4n,2n-1)$, one of the $G_i$ is cyclic of
order $4$ and the other is cyclic of order at least $2n$. To see this,
we can use the description of the action of $\SO(4)$ on $S^3$
described in \cite{Orlik}. Let $R(\theta)$ be the rotation in $\SO(2)$
through an angle $\theta$, and let $M(\theta_1,\theta_2)$ be the
orthogonal sum of $R(\theta_1)$ and $R(\theta_2)$ in $\SO(4)$. Let
$\theta_0$ be $2\pi/(4n)$. Then an element which generates the
$C_{4n}$-action whose quotient is $L(4n,2n-1)$ is
$M(\theta_0,(2n-1)\theta_0)\eq
M(n\theta_0,n\theta_0)M(-(n-1)\theta_0,(n-1)\theta_0)$. The first is
an element of order $4$ in the ``right rotations''. The second is an
element of order~$4n$ or $2n$ in the ``left rotations.'', according
as~$n$ is even or odd.  Therefore $\norm(G^*)$ is $S^1\times S^1$ and
$\isom(L(4n,2n-1))\eq S^1\times S^1$.

\smallskip

\noindent {\sl Case III.} $G\cong D^*_{4m}\times C_n$.

According to \cite{Orlik} we may take $G_1\cong D^*_{4m}$ and
$G_2\cong C_{2n}$.  If $n\eq 1$, then $\norm(G^*)$ is $\Spin(3)$, so
$\isom(M(m,n))\eq\SO(3)$. If $n\geq 3$, then $\norm(G^*)$ and
$\isom(M(m,n))$ are~$S^1$.

\smallskip

\noindent {\sl Case IV.} $G\cong D_{(2^k,m)}\times C_{n_1}$, $k\geq3$.

From the first paragraph on p.~111 of \cite{Orlik} (which applies only
when $k\geq 3$, not for $k\geq 2$ as stated there) we have $G_1\cong
C_{2^kn_1}$ and $G_2\cong D^*_{4m}$. Therefore $\norm(G^*)$ and
$\isom(M(m,n))$ are~$S^1$.

\smallskip

Calculation of $\pi_0(\Diff(M(m,n)))$, implying that
$\pi_0(\Isom(M(m,n)))$ is iso\-morphic to $\pi_0(\Diff(M(m,n)))$, was
done in \cite{A}, \cite{C-S}, and \cite{R2}.  We summarize the
information we have collected so far in the following table, where as
above $k$ and $n_1$ are defined by $4n\eq 2^kn_1$ with $n_1$ odd.

\def\spa{\hbox{\hskip3pt}}

{$$\vbox{\offinterlineskip\halign{
    %%%%%%%%%%%  Preamble   %%%%%%%%%
   \vrule \spa#\spa\hfil&\vrule height 10pt depth 4pt width 0pt \hfil
   \spa#\spa \hfil
&\vrule\hfil \spa#\spa \hfil&\vrule\hfil \spa#\spa \hfil&\vrule\hfil
\spa#\spa \hfil
&\vrule\hfil  \spa#\spa \hfil&\hfil\vrule#\cr
        %%%%  Note   \strut replaced by \vrule height 15pt depth 6pt width 0pt
    %%%%%%%%%%   Chart     %%%%%%%%%%%%
   \noalign{\hrule}
&{\hfil$m,n$ values\hfil}&
       {\hfil$M$\hfil}&
       {\hfil$\pi_1(M)$\hfil}&
       {\hfil$\isom(M)$\hfil}&
       {\hfil$\;\pi_0(\Isom(M))$\hfil}&\cr
\noalign{\hrule}
&&&&&&\cr
\noalign{\vskip-12.5pt}
\noalign{\hrule}

&$m\eq n\eq 1$&$L(4,1)$&$C_4$&$S^1\times \SO(3)$&$C_2$&\cr
\noalign{\hrule}
&$m\eq 1,\,n\neq1$&$L(4n,2n-1)$&$C_{4n}$
&$S^1\times S^1$
&$C_2\times C_2$&\cr
\noalign{\hrule}
&$m\eq 2,\,n\eq 1$&\hbox{quaternionic}&$D^*_8$&$\SO(3)$
&$S_3$&\cr
\noalign{\hrule}
&$m\eq 2,\,n\neq1$&\hbox{prism mfd.}&$D^*_8\times C_n$&
            $S^1$&$S_3\times C_2$&\cr
\noalign{\hrule}
&$m\neq1,2,\,\,n\eq 1$&\hbox{prism mfd.}&$D^*_{4m}$
&$\SO(3)$&$C_2$&\cr
\noalign{\hrule}
&$m\neq 1,2,\,n\neq1$,&\hbox{prism mfd.}&
             $D_{(2^k,m)}\times C_{n_1}$&$S^1$&$C_2\times C_2$&\cr
&$m$ odd\ \ &&&&&\cr
\noalign{\hrule}
&$m\neq 1,2,\,n\neq1$,&\hbox{prism mfd.}&
             $D^*_{4m}\times C_n$&$S^1$&$C_2\times C_2$&\cr
&$m$ even\ \ &&&&&\cr
  \noalign{\hrule}
      }}$$}

\section[Homotopy type of the space of diffeomorphisms]
{Homotopy type of the space of diffeomorphisms}
\label{mainthms}

In the space of smooth imbeddings of the Klein bottle in $M\eq
M(m,n)$, denote by $\imb(K,M)$ the connected component of the
inclusion of the ``standard'' Klein bottle $K_0$, which will be
defined below when we give a more precise description of the Seifert
fiberings we will be using. We will assume that exactly one of $m$ or
$n$ is equal to~$1$.

When $m\eq1$, so that $M$ is a lens space, there is a Seifert fibering
with two exceptional orbits of type $(2,1)$ contained in $K_0$. The
quotient $2$-orbifold is the sphere with two cone points of order $2$.
When $n\eq1$, $M$ is a binary dihedral space and there is a Seifert
fibering which is nonsingular with orbit space $\R\P^2$. In both
cases, $K_0$ is a union of fibers. In the subspace of $\imb(K,M)$
consisting of those imbeddings which take fibers of $K_0$ to fibers of
$M$, let $\imb_f(K,M)$ denote the connected component of the
inclusion. These are called the {\em fiber-preserving} imbeddings.

Our main result shows that parameterized families of imbeddings of $K$
in $M$ can be deformed to families of fiber-preseving imbeddings.

\begin{theorem}{main}{}
If either $m\neq1$ or $n\neq1$, then the inclusion
$\imb_f(K,M)\rightarrow\imb(K,M)$ is a weak homotopy equivalence.
\marginwrite{main}
\end{theorem}

The proof will be given in sections \ref{generic position} and
\ref{parameterization}. From theorem~\ref{main}, we can deduce
the Generalized Smale Conjecture for these classes of 3-manifolds.

\begin{theorem}{smale1}{} For the binary dihedral spaces
$M(m,1)$, $m\geq2$, the inclusion from $\Isom(M(m,1))$ to
$\Diff(M(m,1))$ is a homotopy equivalence, consequently
$\Diff(M(m,1))$ is homotopy equivalent to $\SO(3)\times S_3$ or
$\SO(3)\times C_2$ according as $m\eq2$ or $m>2$.
\marginwrite{smale1}
\end{theorem}

\begin{theorem}{smale2}{} For the lens spaces $M(1,n)\eq L(4n,2n-1)$,
$n\geq2$, the inclusion from $\Isom(M(1,n))$ to $\Diff(M(1,n))$ is a
homotopy equivalence, consequently $\Diff(M(1,n))$ is homotopy
equivalent to $S^1\times S^1\times C_2\times C_2$.
\marginwrite{smale2}
\end{theorem}

Before beginning the proofs, we will need a more precise description
of the Seifert fiberings that are invariant under the
isometries. Assume first that $m\eq1$ so that $M$ is a lens space.  A
generating element of $\pi_1(M)$ was given explicitly in Case~I in
section~\ref{results}, and is a product
$M(n\theta_0,n\theta_0)M(-(n-1)\theta_0,(n-1)\theta_0)$ where
$\theta_0\eq 2\pi/(4n)$. The action of the left rotations contains an
$S^1$-subgroup $S_L$ which contains the cyclic subgroup $C_{2n}$ of
$\pi_1(M)$ generated by the element $M(2\theta_0,-2\theta_0)$;
explicitly, it is the group of elements of the form
$M(\theta,-\theta)$. The orbits of the action of $S_L$ are the fibers
of a Hopf fibering of $S^3$ by geodesic circles. Since the right
rotations commute with the left rotations, the action of the right
rotations is fiber-preserving. The quotient space of the Hopf fibering
is $S^2$ on which the right rotations act via a quotient map $q\colon
R\to
\SO(3)$ described explicitly on p.~105 of \cite{Orlik}. In particular,
$q(M(n\theta_0,n\theta_0))$ is an element of order~2 ($M(\pi,\pi)$ is
the kernel of $q$), which acts on $S^2$ with two fixed points,
corresponding to the two orbits of the $S_L$-action left invariant by
$M(n\theta_0,n\theta_0)$. The quotient of $S^3$ by $C_{2n}$ is
$L(2n,1)$, and the $C_2$-action induced by
$M(\theta_0,(2n-1)\theta_0)$ preserves exactly two orbits which become
the two $(2,1)$ exceptional orbits of $L(4n,2n-1)$. Now let $K_0$ be
the preimage of a great circle of $S^2$ through the two fixed points
of $q(M(n\theta_0,n\theta_0))$. This is a totally geodesic vertical
Klein bottle in $M$. As explained in section~\ref{results},
$\isom(M)\eq S^1\times S^1$ where one $S^1$-factor is the vertical
action on $M$ induced by $S_L$ and the other comes from the
$S^1$-action induced on $M$ by the $S^1$-subgroup $S_R$ that contains
$M(n\theta_0,n\theta_0)$.  Restricted to $K_0$, these isometries give
all the isometric fiber-preserving imbeddings of $K_0$ in $M$: the
vertical reimbeddings are given by the restriction of $S_L$, while the
action of $S_R$ moves $K_0$ through all the Klein bottles that are
preimages of great circles of $S^2$ that pass through the two fixed
points of $q(M(n\theta_0,n\theta_0))$. Summarizing, if we denote the
isometric fiber-preserving imbeddings of $K_0$ by $\isom_f(K,M)$, we
have shown that the restriction map $\isom(M)\to\isom_f(K,M)$ is a
homeomorphism.

Suppose now that $n\eq 1$. From section~\ref{results}, $\pi_1(M)\cong
D^*_{4m}$ and from Case~III of section~\ref{results} we may assume
that $D^*_{4m}$ is a subgroup of the right rotations. The Seifert
fibering invariant under $\isom(M)$ is obtained as follows (see
pp.~112-113 of~\cite{Orlik}). There is an $S^1$-subgroup $S$ of the
group of right rotations which contains the index~2 subgroup
$C^*_{2m}$ of $D^*_{4m}$. There is an order~4 right rotation $\delta$
which conjugates each element of $S$ to its inverse, and $D^*_{4m}$ is
generated by $C^*_{2m}$ and $\delta$. The orbits of $S$ are preserved
by $\delta$ and determine the fibering of $M$. Now let $p\colon S^3\to
S^2$ be the Hopf map whose point preimages are the orbits of $S$. On
$S^2$, $\delta$ induces the antipodal map. Fix one of the invariant
circles of $\delta$. Its image under $p$ is a great circle $C$, and we
let $K_0$ be the image of $p^{-1}(C)$ in $M$.  The induced fibering on
$K_0$ is the nonsingular one by meridinal fibers, and for this metric
on $K_0$, $\isom_f(K_0)$ is $S^1$ in which the order~2 element takes
each fiber to itself by the monodromy of $K_0$ regarded as a
$S^1$-fibering over $p(C)$. As seen in section~\ref{results}, the
isometries of $M$ are induced by the left rotations. In particular,
the circle subgroup of the group of left rotations that leaves
$p^{-1}(C)$ invariant restricts to $\isom_f(K_0)$ on~$K_0$. As in the
case $m\eq 1$, the restriction determines a homeomorphism $\isom(M)\to
isom_f(K,M)$.

\begin{proof}{Theorems \ref{smale1} and \ref{smale2} assuming Theorem
\ref{main}} Since $\Diff(M)$ has the homotopy type of a CW-complex
\cite{P}, it is enough to prove that the inclusion is a weak homotopy
equivalence. As mentioned in section~\ref{results}, the inclusion is
known to induce a bijection on path components, so we will restrict
attention to the connected components of the identity homeomorphism.

From corollary~8.7 of \cite{Kalliongis-McCullough}, restriction of
diffeomorphisms to imbeddings defines a fibration
$$\Diff_f(M\rel K_0)\cap \diff_f(M)\to \diff_f(M)\to \imb_f(K,M)\ .$$

\noindent Since any diffeomorphism in this fiber is
orientation-preserving, it cannot interchange the sides of
$K_0$. Therefore the fiber may be identified with a subspace
consisting of path components of $\Diff_f(S^1\times D^2\rel S^1\times
\partial D^2)$. By theorem~5.2 of~\cite{Kalliongis-McCullough}, there
is a fibration
$$\Diff_v(S^1\times D^2\rel S^1\times D^2)\to \Diff_f(S^1\times
D^2\rel S^1\times D^2)\to \Diff(D^2\rel\partial D^2)\ ,$$

\noindent whose fiber is the group of {\em vertical} diffeomorphisms
that take each fiber to itself. The base is contractible
by~\cite{Smale}. The fiber is contractible, this is seen by lifting
diffeomorphisms to the infinite cyclic cover $\R\times D^2$ and
canonically and equivariantly deforming the lifts to preserve
$\set{0}\times D^2$, and then to be the identity. We conclude that
$\Diff_f(S^1\times D^2\rel S^1\times \partial D^2)$ and hence also
$\Diff_f(M\rel K_0)$ are contractible.
Therefore our fibration from above becomes
$$\diff_f(M\rel K_0)\to \diff_f(M)\to \imb_f(K,M)$$

\noindent with contractible fiber. Similarly there is a fibration
$$\diff(M\rel K_0)\to \diff(M)\to \imb(K,M)\ .$$

\noindent The fact that it is a fibration  comes from \cite{Palais} and the
contractibility of the fiber uses~\cite{H2}. We can now fit these into
a diagram
$$\vbox{\halign {\hfil$#$\hfil&\hfil$#$\hfil&
        \hfil$#$\hfil&\hfil$#$\hfil&\hfil$#$\hfil&\hfil$#$\hfil&
        \hfil$#$\hfil&\hfil$#$\hfil&\hfil$#$\hfil\cr
        & & \diff_f(M\rel K_0)
         & \mapright{} & \diff_f(M)
         & \mapright{} & \imb_f(K,M)
         &  &  \cr
         &  & \mapdown{}
         &  & \mapdown{}
         &  & \mapdown{}
         &  &  \cr
        &  & \diff(M\rel K_0)
         & \mapright{} & \diff(M)
         & \mapright{} & \imb(K,M)\ .
         &  &  \cr}}$$

\noindent The vertical maps are inclusions. By theorem~\ref{main},
the right hand vertical arrow is a weak homotopy equivalence. Since
the fibers are both contractible, it follows that $\diff_f(M)\to
\imb_f(K,M)$, $\diff(M)\to \imb(K,M)$, and $\diff_f(M)\to \diff(M)$
are weak homotopy equivalences.

Let ${\cal O}_K$ be the image of $K_0$ in the quotient orbifold ${\cal
O}$ of the fibering on $M$. When $m\eq1$, ${\cal O}_K$ is a
silvered interval imbedded as half of a great circle connecting the
two order~2 cone points of ${\cal O}$. When $n\eq1$, ${\cal O}_K$ is an
$\R\P^1$ which is the image of a great circle of $S^2$ in ${\cal O}\eq
\R\P^2$. Each element of $\isom_f(K,M)$ projects to an isometric
imbedding of orbifolds of ${\cal O}_K$ in~${\cal O}$.

Denote by $\imb({\cal O}_K,{\cal O})$ the connected component of the inclusion
in the space of orbifold imbeddings, and let a subscript $v$ as in
$\Diff_v(K_0)$ indicate the vertical maps that take each fiber to itself.
We have a fibration of groups
$$\Isom_v(K_0)\cap \isom_f(K,M)\to \isom_f(K,M)\to \isom({\cal
O}_K,{\cal O})\ .$$

\noindent When $m\eq 1$, this sequence is readily seen to be
$S^1\to S^1\times S^1\to S^1$. When $n\eq 1$, $\isom({\cal O}_K,{\cal
O})$ can be identified with the unit tangent space of $\R\P^2$, and
$\Isom_v(K_0)\cap\isom_f(K,M)\eq C_2$ generated by applying the
monodromy (of the $S^1$-bundle $K_0\to p(C)$ described above) in each
fiber. In this case, the sequence is $C_2\to \SO(3)\to T_1(\R\P^2)$
(topologically this is $C_2\to \R\P^3\to L(4,1)$). Using inclusions as
the vertical maps, we have a diagram of fibrations
$$\vbox{\halign {\hfil$#$\hfil&
        \hfil$#$\hfil&\hfil$#$\hfil&\hfil$#$\hfil&\hfil$#$\hfil&
        \hfil$#$\hfil&\hfil$#$\hfil&\hfil$#$\hfil&\hfil$#$\hfil\cr
        & & \Isom_v(K_0)\cap \isom_f(K,M)
         & \mapright{} & \isom_f(K,M)
         & \mapright{} & \isom({\cal O}_K,{\cal O})
         &  &  \cr
         &  & \mapdown{}
         &  & \mapdown{}
         &  & \mapdown{}
         &  &  \cr
        &  & \Diff_v(K_0)\cap \imb_f(K,M)
         & \mapright{} & \imb_f(K,M)
         & \mapright{} & \imb({\cal O}_K,{\cal O})
         &  &  \cr}}$$

\noindent where the bottom sequence is a fibration by theorem~8.9 of
\cite{Kalliongis-McCullough}.

Suppose first that $n\eq1$. Using theorem 4 of \cite{G}, the
right-hand vertical arrow of the diagram is a homotopy
equivalence. The two components of $\Diff_v(K_0)$ are contractible,
each containing a unique isometry, and thus the left-hand vertical
arrow is a homotopy equivalence. It follows that the middle vertical
arrow is a weak homotopy equivalence. When $m\eq1$, $\Isom_v(K_0)$ and
$\Diff_v(K_0)$ have two components, but elements of one component
reverse the direction of the fiber so are not contained in
$\imb_f(K,M)$. The identity component $\isom_v(K_0)$ is homeomorphic
to $S^1$, and the inclusion $\isom_v(K_0)\to \diff_v(K_0)$ is a
homotopy equivalence, although the full details of this are
lengthy. Again the middle arrow is a weak homotopy equivalence.

We have seen that $\isom(M)\to isom_f(K,M)$ is a
homeomorphism. (This fails when $(m,n)\eq(1,1)$, since then
$\isom(M)$ does not preserve any Seifert fibering of $M$.)
We now have a diagram of inclusions

$$\vbox{\halign {\hfil$#$\hfil&\hfil$#$\hfil&
        \hfil$#$\hfil\cr
         \isom(M)
         & \mapright{} & \isom_f(K,M)
         \cr
        \mapdown{}
         & & \mapdown{}
          \cr
         \diff_f(M)
         & \mapright{} & \imb_f(K,M)
          \cr
         \mapdown{}
         & & \mapdown{}  \cr
         \diff(M)
         & \mapright{} & \imb(K,M)
          \cr}}$$

\noindent in which all arrows except the one from $\isom(M)$ to $\diff_f(M)$
have been shown to be weak homotopy equivalences; it follows that it
and the composite $\isom(M) \to \diff(M)$ are weak homotopy equivalences
as well.
\end{proof}

\section [Generic Position Configurations] {Generic Position
Configurations}
\label{generic position}\marginwrite{generic position}

A smoothly imbedded (connected) 2-manifold $T$ in a closed 3-manifold
$M$ has either a product neighborhood $T\times [-1,1]$, or a 2-fold
covering from $\widetilde T\times [-1,1]$ to a tubular neighborhood of
$T$. In the former case, let $T_u$ denote $T\times \{u\}$, and in the
latter let $T_u$ denote the image of $\widetilde T\times\{u\}$. We
call the $T_u$ {\it horizontal levels} of the neighborhood of $T$.

Let $S$ and $T$ be smoothly imbedded closed surfaces in a closed
3-manifold $M$.  A point $x$ in $S\cap T$ is called a {\it regular}
point if $S$ is transverse to $T$ at $x$, otherwise it is a {\it
singular} point. Following section~5 of \cite{I2}, we call $x$ a {\it
singular point of finite multiplicity} if $S\cap T$ meets a small
neighborhood $U$ of $x$ in a finite even number of smooth arcs running
from $x$ to $\partial U$, transversely except at $x$ (cf.~Fig.\ 3, p.\
1653 of \cite{I2}). Then, either $S\cap T\cap U\eq \{x\}$ or $x$ is a
saddle tangency of $S$ and~$T$.

We say that the surfaces are {\it in generic position} if all singular
points of intersection are of finite multiplicity. A parameterized
family in $\Imb(S,M)$ is said to be in generic position relative to
$T$ if each of the imbeddings in the family has this property.

\begin{proposition}{perturb}{} Suppose that $F\colon  D^k\to\hbox{\it
Imb}(S,M)$ is a parameterized family of imbeddings. Assume either that
$F(t)(S)$ is in generic position relative to $T$ for all $t$ in
$\partial D^k$, or that $F(t)(S)\eq T$ for all $t$ in $\partial
D^k$. Then $F$ is homotopic relative to $\partial D^k$ to a map
$G\colon D^k\to\hbox{\it Imb}(S,M)$ so that $G\eq F$ on $\partial D^k$
and $G(t)(S)$ is in generic position with respect to $T$ for all $t\in
\hbox{\it int}(D^k)$. Moreover, for each $t\in \hbox{\it int}(D^k)$
there exists $u_0>0$ so that $G(t)(S)$ is transverse to $T_u$ for all
$0<u\leq u_0$, where $T_u$ are horizontal levels in a tubular
neighborhood of $T$.
\marginwrite{perturb}
\end{proposition}

\noindent For a discussion of this proposition, we refer the reader to
lemma (5.2) and remark (5.3) of \cite{I2}. The map $G$ may be chosen
arbitrarily close to $F$, although we will not need to do so.

Suppose now that $L_0$ is a 1-sided surface in $M$, and as above let
$L_u$ denote the horizontal levels of a tubular neighborhood of $L_0$.
A piecewise-linearly imbedded surface $S$ in $M$ is said to be {\it
flattened} (with respect to $L_0$ and the choice of the  $L_u$) if it
satisfies the following conditions.

\begin{enumerate}
\item[1.] There is a 4-valent graph $\Gamma$ (possibly with components
          which are circles) contained in $L_0$ such
          that $S\cap L_0$ consists of the closures
          of some of the connected components of $L_0-\Gamma$.
\item[2.] Each point $p$ in the interior of an
             edge of $\Gamma$
           has a neighborhood $U$ for which the quadruple
          $(U,U\cap L_0,U\cap S,p)$ is PL homeomorphic to the configuration
          $(\R^3,\{(x,y,z)\;\vert\;z\!\eq \!0\},
          \{(x,y,z)\;\vert\;\hbox{{\rm either\ }}z\!\eq \!0
          \hbox{{\rm\ and\ }}y\geq0,\hbox{{\rm\ or\ }}
          y\!\eq \!0\hbox{{\rm\ and\ }}z\geq0\},\{0\})\,$
          (see Figure $1$(a)).
\item[3.] Each vertex $v$
          of $\Gamma$ has a neighborhood $U$ for which the quadruple
          $(U,U\cap L_0, U\cap S,v)$ is PL homeomorphic to the configuration
          $(\R^3,\{(x,y,z)\;\vert\;z\eq 0\},
          \{(x,y,z)\;\vert\;\hbox{{\rm either\ }}z\!\eq \!0
          \hbox{{\rm\ and\ }}xy\leq0,\hbox{{\rm\ or\ }}
          y\!\eq \!0\hbox{{\rm\ and\ }}z\geq0,\hbox{{\rm\ or\ }}x\!\eq \!0
          \hbox{{\rm\ and\ }}z\leq0\},\{0\})\,$ (see Figure $1$(b)).
\end{enumerate}

\begin{lemma}{flatten}{} Let $N$ be a 3-manifold containing a smoothly
imbedded surface $L_0$, and let $S_1$ be a smoothly imbedded surface
in $N$ which meets $L_0$ in generic position and meets each $L_u$
transversely, for $0<u\leq u_0$. Then given $\epsilon>0$ there is a PL
isotopy $S_t$ from $S_1$ to a PL imbedded surface $S_0$ such that
\begin{enumerate}
\item[{\rm(1)}] Each $S_t$ is within  distance $\epsilon$ of the inclusion
         $S_1$.
\item[{\rm(2)}] Each $S_t$ is transverse to $L_u$ for $0<u\leq u_0$.
\item[{\rm(3)}] $S_0$ is flattened.
\end{enumerate}
\marginwrite{flatten}
\end{lemma}

\begin{proof}{} The isotopy will move points monotonically with respect
to $u$ levels. We first describe it near a singular point $x$ of
$S_1\cap L_0$. In a neighborhood $U$ of $x,\;S_1\cap L_0$ consists of
$x$ together with a (possibly empty) collection of arcs
$\alpha_1,\alpha_2,\ldots,\alpha_{2n}$ emanating from $x$. With
respect to some fixed Riemannian metric for $N$, there is a
neighborhood of $x$ for which the angle of intersection of $S_1$ with
$L_0$ is small; the isotopy moves points only within an $\epsilon$
neighborhood of the $\alpha_i$ and decreases these angles to $0$
everywhere in a neighborhood of $x$ (or pushes a $2$-disc neighborhood
of $x$ in $S_1$ down to a $2$-disc neighborhood of $x$ in $L_0$, if
there are no arcs). At the end of the initial isotopy, say for $1\geq
t\geq 1/2$, there is a neighborhood $U$ of $x$ for which $S_{1/2}\cap
L_0\cap U$ is a regular neighborhood in $L_0$ of
$\cup_{i=1}^{2n}\alpha_i$. These isotopies may be performed
simultaneously near each singular point of intersection. The remainder
of the isotopy will take place in an $\epsilon$ neighborhood of the
original (open) edges of $S_1\cap L_0$. At the end of this isotopy,
the intersection will be locally a regular neighborhood of the
original edges, except that on some of the edges it might be necessary
to introduce a point where the configuration is as in Figure 1(b)---
this is necessary only when the flattenings at the singular points at
the ends of the edge are in opposite senses. Again, these remaining
isotopies may be performed in disjoint neighborhoods of the original
edges.
\end{proof}

We call an isotopy as in lemma~\ref{flatten} a {\it flattening
isotopy.} By property (2), the collection of intersection circles in
$L_u$ for $0<u\leq u_0$ is changed only by isotopy in $L_u$. After
flattening, each of these circles projects along $S_0$ to an immersed
circle lying in $\Gamma$, having a transverse self-intersection at
each of its double points (which can occur only at vertices of
$\Gamma$.)

Now we specialize to the standard Klein bottle $K_0\subseteq M$ and a
parameterized family $F\colon D^k\to \imb(K,M)$. We denote the
imbedding $F(t)$ by $F_t$ and its image by $K_t$. Assume for all
$t\in\partial D^k$ either that $K_t\eq K_0$, or that $K_t$ is
transverse to $K_0$.  By proposition~\ref{perturb}, we may deform $F$
to be in generic position relative to $K_0$ over $\hbox{\it
int}(D^k)$.

\begin{theorem}{circles}{} Suppose that $M\eq M(m,n)$ with one but not
both of $m$ or $n$ equal to 1, and let $K_t$ be a parameterized family
in generic position for each $t\in\hbox{\it int}(D^k)$. Then for each
$t\in\hbox{\it int}(D^k)$, there exists $u>0$ so that $K_t$ is
transverse to $T_u$ and each circle of $K_t\cap T_u$ is either
inessential in $T_u$, or represents $a$ or $b^2$ in $\pi_1(T_u)$.
\marginwrite{circles}
\end{theorem}

It follows immediately that no circle of $K_t\cap T_u$ is homotopic in
$T_u$ to the meridian. Moreover, it follows that no intersection
circle is homotopic in $T_u$ to a longitude of $R_u$ which is not
homotopic in $T_u$ to a fiber of the Seifert fibering. For when $n\eq
1$, the longitudes are $a(a^mb^2)^k$, where $k$ is an arbitrary
integer, and when $m\eq 1$ they are $b^2(ab^{2n})^k$. (In $L(4,1)$,
however, an $a$ circle is a longitude of $R_u$ which is not homotopic
to a fiber of the fibering with two $(2,1)$ orbits, while a $b^2$
circle is a longitude not homotopic to a fiber of the nonsingular
fibering.)

The proof will produce $u$ so that $K_t$ is transverse
to $T_s$ for $0<s\leq u$, but we will not need this property.

\begin{proof}{} Suppose first that the intersection $K_t\cap K_0$ is
transverse. Since $K_t$ must meet every nearby level $T_u$
transversely, it intersects $P_u$ in M\"obius bands and annuli.
Consequently the projection of $T_u$ onto $K_0$ maps circles of
intersection of $K_t\cap T_u$ onto circles of $K_t\cap K_0$ either
homeomorphically or by two-fold coverings. Only inessential and $a$
and $b^2$ circles can be preimages of imbedded circles in~$K_0$.

Suppose now that $K_t$ intersects $K_0$ in some singular points. By
proposition~\ref{perturb}, $K_t$ is transverse to $T_u$ for all $u\leq
u_0$ for some $u_0$. By lemma~\ref{flatten}, we can flatten $K_t$ near
$K_0$, without changing either the transversality of $K_t$ and $T_u$ or
the homotopy classes of the loops in $K_t\cap T_u$. Then, $K_t\cap K_0$
consists of a valence 4 graph $\G$, which is the image of the
collection of disjoint simple closed curves $K_t\cap T_u$ under a
2-fold covering projection, together with some of the complementary
regions of $\G$ in $K_0$, which we will call the {\it faces.} Each edge
of $\G$ lies in exactly one face. We may choose the $I$-fibering so
that $K_t\cap P_u$ lies in the union of $K_t\cap K_0$ and the
$I$-fibers that meet $\G$.

Suppose for contradiction that one of the circles in $K_t\cap T_u$
represents $a^kb^{2\ell}$ with $k\ell\neq0$. Since $K_t$ is
geometrically incompressible (if not, then $M$ would contain an
imbedded projective plane), there is an isotopy of $K_t$ which
eliminates the circles of $K_t\cap T_u$ that are inessential in $T_u$,
without altering the remaining circles or destroying the flattened
position of $K_t\cap P_u$. So we may assume that $K_t\cap T_u$
consists of disjoint circles each representing $a^kb^{2\ell}$. Since
$K_t$ is isotopic to $K_0$, each loop in $T_u$ has even algebraic
intersection number with $K_t\cap T_u$, so there is an even number of
these circles; denote them by $A_1$, $A_2,\ldots\,$, $A_{2r}$.

The vertices of $\G$ are the images of the intersections of $\cup A_i$
with $\cup\t(A_i)$, where the involution $\t$ is the covering
transformation for the covering map from $T_u$ to $K$ that determines
the $I$-fibering of $P_u$. Now $\cup A_i$ and $\cup \t(A_i)$ meet
transversely; the number of intersections is at least
$$\abs{\left(\cup A_i\right)\,\cdot\,\left(\cup \t(A_i)\right)}
   = \abs{(2r\,a^kb^{2\ell})\,\cdot\,(2r\,a^kb^{-2\ell})}
   = 4r^2\,\abs{2k\ell}\ .$$

\noindent Since each vertex of $\G$ is covered by two intersections,
$\G$ has at least $4r^2\abs{k\ell}$ vertices.

Notice that each edge of $\G$ runs between two distinct vertices,
since the $A_i$ are disjoint and imbedded and do not cover imbedded
loops in $K$. Moreover, each face contains an even number of edges,
since it can be lifted to $T_u$ with successive edges lying
alternately in $\cup A_i$ and $\cup \t(A_i)$. In particular, no face
is a 1-gon, or is a 2-gon with its vertices identified. Therefore each
face that is a 2-gon can be eliminated by an isotopy, yielding a new
$K_t$ in flattened position (see Fig.~2). So we may assume that each
face contains at least 4 vertices. Finally, observe that the Euler
characteristic of $K_t\cap P_u$ is at least $-2r$, since
$\chi(K_t)\eq0$ and $K_t\cap P_u$ has exactly $2r$ boundary
components. Letting $V$, $E$, and $F$ denote the number of vertices,
edges, and faces of $K_t\cap K_0$, we have $E\eq2V$ and $F\leq V/2$
(since each edge lies in exactly one face and each face has at least
$4$ edges). Therefore $-2r\leq \chi(K_t\cap P_u)\eq\chi(K_t\cap
K)\leq-V/2$. Since $V\geq 4r^2\abs{k\ell}$, it follows that
$r\abs{k\ell}\leq1$, forcing $r\eq\abs{k\ell}\eq1$, $\chi(K_t\cap
K)\eq -2$, $V\eq 4$, and $F\eq2$.  That is, $K_t\cap K$ consists of
two faces, each a $4$-gon, meeting at their four vertices. Moreover, $\G$
is the image of two imbedded circles $ab^2$ or $ab^{-2}$ circles each
projecting to a loop with one self-intersection. This forces the two
faces of $K_t\cap K_0$ to meet each other as shown in Fig.~3. It
follows that $K_t\cap P_u$ is a twice-punctured Klein bottle (shading
the complementary faces would yield a twice-punctured torus). But then
$ab^2$ or $ab^{-2}$ bounds a disc in $R_u$, contradicting the fact
that $ab^2$ and $ab^{-2}$ are nontrivial in $\pi_1(M(m,n))$ (since
$(m,n)\neq(1,1)$).
\end{proof}

\section[Parameterization]{Parameterization}
\label{parameterization}

We now complete the proof of theorem~\ref{main}. Since $\imb(K,M)$
and $\imb_f(K,M)$ are connected, we have
$\pi_0(\imb(K,M),\imb_f(K,M))\eq0$. To prove that the higher relative
homotopy groups vanish, we begin with a parameterized family, which we
may take to be a smooth map $D^k\rightarrow \imb(K,M)$, where
$k\geq1$, which takes all points of $\partial D^k$ to the standard
inclusion. By abuse of notation, we confuse the imbedding
corresponding to the point $t\in D^k$ with its image, denoting both by
$K_t$. By proposition~\ref{perturb}, we may assume that each $K_t$ is
in general position with $K_0$. By theorem~\ref{circles} and the
remark following its statement, there is for each $t$ a value $u>0$ so
that

\begin{enumerate}
\item[(1)] $K_t$ is transverse to $T_u$.
\item[(2)] No intersection circle of $K_t$ with $T_u$
           is a meridian, or a longitude not
           homotopic in $T_u$ to a fiber of the Seifert fibering.
\end{enumerate}

\noindent Each intersection circle that bounds a (necessarily unique)
2-disc in $T_u$ also bounds a unique 2-disc in $K_t$, since $K_t$ is
geometrically incompressible, and if either of the two discs is
innermost among all such discs, then their union bounds a unique
3-ball in $M$. Only very routine modifications are needed to the
procedure of Hatcher described in \cite{H3} to deform the family,
keeping it fixed on $\partial D^k$,  so that for each $t\in D^k$,
there is a value $u>0$ so that in addition to (1) and (2) we have

\begin{enumerate}
\item[(2$^{\prime}$)] No intersection circle of $K_t$
      with $T_u$ is inessential in $T_u$.
\end{enumerate}

\noindent In fact, the argument is somewhat simpler, since there is
only a unique 3-ball across which the 2-discs can be pushed to
eliminate the inessential circles; moreover since these 3-balls cannot
contain essential loops of $T_u$ in their boundaries, the deformations
can be chosen so as not to affect the intersections which are
essential in the $T_u$. As in \cite{H3}, it is necessary to pass to
new levels, but by choosing these very close to previously chosen
levels we can ensure that no new kinds of intersection circles arise.

Since $K_t$ is incompressible, it follows that $K_t\cap R_u$ consists
of annuli (M\"obius bands cannot occur because orientation-reversing
loops in $K_t$ are dual to $[K_0]\in H_2(M;\Z/2))$.

Consider the annuli $K_t\cap R_u$ whose boundary circles are not
longitudes. Each such annulus is parallel in $R_u$ to a uniquely
determined annulus in $T_u$. We again use the procedure of \cite{H3}
to pull these annuli out of the $R_u$. That is, deform the family,
keeping it fixed on $\partial D^k$,  so that for each $t\in D^k$,
there is a value $u>0$ so that

\begin{enumerate}
\item[(1)] $K_t$ is transverse to $T_u$.
\item[(2$^{\prime\prime}$)] Every intersection circle of $K_t$ with $T_u$
      is homotopic in $T_u$ to a fiber of the Seifert fibering.
\end{enumerate}

\noindent The adaptation of the argument is routine; the annuli play
the role of the regions called $D_M(c_i)$ in Hatcher's paper, and the
cross-sectional picture of the regions between the annuli and the $T_u$
is exactly as in the figure on p.~429 of \cite{H3}.
Again, it is necessary to pass to new levels, but
no new kinds of intersection circles need arise.

To complete the argument, we require two technical lemmas.

\begin{lemma}{fiberpreserving1}{} Let $T$ be a torus with a fixed
$S^1$-fibering, and let $C_n\eq\cup_{i=1}^nS_i$ be a union of $n$
distinct fibers.  Then $\imb_f(C_n,T)\to\imb(C_n,T)$ is a weak
homotopy equivalence.
\marginwrite{fiberpreserving1}
\end{lemma}

\begin{proof}{} Fix a basepoint $s_0$ in $S_n$.
Consider the diagram
$$\vbox{\halign {\hfil$#$\hfil&\hfil$#$\hfil&\hfil$#$\hfil&\hfil$#$\hfil&%
        \hfil$#$\hfil&\hfil$#$\hfil&\hfil$#$\hfil\cr
       & & \imb_f(S_n,T\rel s_0)
         &\to
         & \imb_f(S_n,T)
         &\to
         & \imb(s_0,T)
         \cr
         &  & \mapdown{}
         &  & \mapdown{}
         &  & \mapdown{=}
         \cr
      &  & \imb(S_n,T\rel s_0)
         &\to
         & \imb(S_n,T)
         &\to
         & \imb(s_0,T)
         \cr
}}$$

\noindent The top row is a fibration by corollary~9.6
of~\cite{Kalliongis-McCullough}, and the bottom row is a fibration
by~\cite{Palais}. The fiber of the
first row is homeomorphic to $\Diff_+(S_n\rel s_0)$, the group of
orientation-preserving diffeomorphisms, which is contractible. The
fiber of the second row is contractible using~\cite{G}.
Therefore the middle vertical arrow is a weak homotopy
equivalence. For $n\eq 1$, this completes the proof. Inductively, let
$A$ be the annulus that results from cutting $T$
along $S_n$ and consider a similar diagram, where $\imb(S_{n-1},A)$
denotes the imbeddings with image in the interior of~$A$, and so on.
$$\vbox{\halign {\hfil$#$\hfil&\hfil$#$\hfil&\hfil$#$\hfil&\hfil$#$\hfil&
        \hfil$#$\hfil\cr
          \imb_f(S_{n-1}, A\rel s_0)
         &\to
         & \imb_f(S_{n-1},A)
         &\to
         & \imb(s_0,\hbox{\it int}(A))
         \cr
          \mapdown{}
         & & \mapdown{}
         & & \mapdown{=}
         \cr
           \imb(S_{n-1}, A\rel s_0)
         &\to
         & \imb(S_{n-1},A)
         &\to
         & \imb(s_0,\hbox{\it int}(A))
         \cr
}}$$

\noindent The fibers are contractible, so the middle vertical arrow is
a weak homotopy equivalence.  Now we examine another diagram.
$$\vbox{\halign {\hfil$#$\hfil&\hfil$#$\hfil&\hfil$#$\hfil&
        \hfil$#$\hfil&\hfil$#$\hfil\cr
         \imb_f(C_{n-1},A\rel S_{n-1})
         & \to
         & \imb_f(C_{n-1},A)
         & \to
         & \imb_f(S_{n-1},A)
         \cr
         \mapdown{}
         &  & \mapdown{}
         &  & \mapdown{}
         \cr
           \imb(C_{n-1},A\rel S_{n-1})
         & \to
         & \imb(C_{n-1},A)
         & \to
         & \imb(S_{n-1},A)
         \cr
}}$$

\noindent The first row is a fibration by corollary~6.5
of~\cite{Kalliongis-McCullough} and the second is a fibration
by~\cite{Palais}.  The right vertical arrow was shown to be a weak
homotopy equivalence by the previous diagram, and the left one is a
weak homotopy equivalence by induction, so the middle one is also. The
proof is now completed by the diagram
$$\vbox{\halign {\hfil$#$\hfil&\hfil$#$\hfil&
        \hfil$#$\hfil&\hfil$#$\hfil&\hfil$#$\hfil\cr
          \imb_f(C_n,T\rel S_n)
         & \to
         & \imb_f(C_n,T)
         & \to
         & \imb_f(S_n,T)
         \cr
         \mapdown{}
         &  & \mapdown{}
         &  & \mapdown{}
         \cr
           \imb(C_n,T\rel S_n)
         & \to
         & \imb(C_n,T)
         & \to
         & \imb(S_n,T)\ .
         \cr
}}$$
\end{proof}

\begin{lemma}{fiberpreserving2}{}
Let $\Sigma$ be a compact 3-manifold with nonempty boundary and having
a fixed Seifert fibering.  Let $F$ be a compact 2-manifold properly
imbedded in $\Sigma$, such that $F$ is a union of fibers. Let
$\imb_{\partial f}(F,\Sigma)$ be the connected component of the
inclusion in the space of (proper) imbeddings for which the image of
$\del F$ is a union of fibers. Then $\imb_f(F,\Sigma)\to
\imb_{\partial f}(F,\Sigma)$ is a weak homotopy equivalence.
\marginwrite{fiberpreserving2}
\end{lemma}

To prove lemma~\ref{fiberpreserving2}, we need a preliminary result.

\begin{lemma}{sublemma for fiberpreserving2}{} The following maps
induced by restriction are fibrations.
\begin{enumerate}
\item[{\rm(i)}] $\imb(F,\Sigma)\to \imb(\partial F,\partial \Sigma)$
\item[{\rm(ii)}] $\imb_{\partial f}(F,\Sigma)\to \imb_f(\partial
F,\partial \Sigma)$
\item[{\rm(iii)}] $\imb_f(F,\Sigma)\to \imb_f(\partial F,\partial \Sigma)$.
\end{enumerate}
\end{lemma}

\begin{proof}{} Part (ii) follows from part (i) since
$\imb_{\partial f}(F,\Sigma)$ is the preimage of $\imb_f(\partial
F,\partial\Sigma)$ under the fibration of part (i). Parts~(i)
and~(iii) are cases of corollaries~9.3 and~9.4
of~\cite{Kalliongis-McCullough}.
\end{proof}

\begin{proof}{\ref{fiberpreserving2}}
First we use the following fibration from theorem~8.3
of~\cite{Kalliongis-McCullough},
$$\Diff_v(\Sigma\rel\partial\Sigma)\cap
\diff_f(\Sigma\rel\partial\Sigma)
\to\diff_f(\Sigma\rel\partial\Sigma)
\to \diff({\cal O}\rel\partial {\cal O})$$

\noindent where ${\cal O}$ is the quotient orbifold of $\Sigma$ and
as usual $\Diff_v$ indicates the diffeomorphisms that take each fiber
to itself. The orbifold diffeomorphism group of ${\cal O}$ is homotopy
equivalent to a subspace consisting of path components of the
diffeomorphism group of the 2-manifold $B$ obtained by removing the
cone points from~${\cal O}$. Since $\partial B$ is nonempty,
$\diff(B\rel\partial B)$ and therefore $\diff({\cal O}\rel\partial
{\cal O})$ are contractible. Moreover, since $\pi_1(\diff({\cal
O}\rel\partial {\cal O}))$ is trivial, the homotopy exact sequence of
the fibration shows that $\Diff_v(\Sigma\rel\partial\Sigma)\cap
\diff_f(\Sigma\rel\partial\Sigma)$ is connected so equals
$\diff_v(\Sigma\rel\partial\Sigma)$. It is not difficult to see that
each component of $\Diff_v(\Sigma\rel\partial\Sigma)$ is contractible
(see lemma~10.4 of \cite{Kalliongis-McCullough} for a similar argument),
so we conclude that $\diff_f(\Sigma\rel\partial\Sigma)$ is weakly
contractible.

Next, consider the diagram
$$\vbox{\halign {\hfil$#$\hfil&\hfil$#$\hfil&
        \hfil$#$\hfil&\hfil$#$\hfil&\hfil$#$\hfil\cr
         \Diff_f(\Sigma\rel F\cup \partial \Sigma)\cap
               \diff_f(\Sigma\rel\partial\Sigma)
         & \to
         & \diff_f(\Sigma\rel\partial\Sigma)
         & \to
         & \imb_f(F,\Sigma\rel\partial F)
         \cr
         \mapdown{}
         &  & \mapdown{}
         &  & \mapdown{}
         \cr
         \Diff(\Sigma\rel F\cup\partial\Sigma)\cap
               \diff(\Sigma\rel\partial\Sigma)
         & \to
         & \diff(\Sigma\rel\partial\Sigma)
         & \to
         & \imb(F,\Sigma\rel\partial F)
         \cr
}}$$

\noindent where the rows are fibrations by corollaries~8.7 and~3.6
of~\cite{Kalliongis-McCullough}. From above, the components of
$\Diff_f(\Sigma\rel\partial\Sigma)$ and (by cutting along $F$) the
components of $\Diff_f(\Sigma\rel F\cup\partial\Sigma)$ are weakly
contractible. By~\cite{H2}, the components of
$\Diff(\Sigma\rel\partial\Sigma)$ and $\Diff(\Sigma\rel
F\cup\partial\Sigma)$ are weakly contractible. Therefore
to show that $\imb_f(F,\Sigma\rel\partial
F)\to\imb(F,\Sigma\rel\partial F)$ is a weak homotopy equivalence it
is sufficient to show that $\pi_0(\Diff_f(\Sigma\rel F\cup \partial
\Sigma)\cap \diff_f(\Sigma\rel\partial\Sigma))\to
\pi_0(\Diff(\Sigma\rel F\cup\partial\Sigma)\cap
\diff(\Sigma\rel\partial\Sigma))$ is bijective. It is surjective
because every diffeomorphism of a Seifert-fibered 3-manifold which is
fiber-preserving on the (non-empty) boundary is isotopic relative to
the boundary to a fiber-preserving diffeomorphism (lemma~VI.19 of
\cite{Jaco}). It is injective because fiber-preserving diffeomorphisms
that are isotopic are isotopic through fiber-preserving
diffeomorphisms (see~\cite{Waldhausen}).

The proof is completed by the following diagram in which the rows
are fibrations by parts (iii) and (ii) of lemma~\ref{sublemma for
fiberpreserving2}, and we have verified that the left vertical arrow
is a weak homotopy equivalence.
$$\vbox{\halign {\hfil$#$\hfil&\hfil$#$\hfil&
        \hfil$#$\hfil&\hfil$#$\hfil&\hfil$#$\hfil\cr
         \imb_f(F,\Sigma \rel\partial F)
         & \to
         & \imb_f(F,\Sigma)
         & \to
         & \imb_f(\partial F,\partial \Sigma)
         \cr
         \mapdown{}
         &  & \mapdown{}
         &  & \mapdown{=}
         \cr
         \imb(F,\Sigma\rel\partial F)
         & \to
         & \imb_{\partial f}(F,\Sigma)
         & \to
         & \imb_f(\partial F,\partial \Sigma)\ .
         \cr
}}$$
\end{proof}

We can now complete the proof of theorem~\ref{main} by deforming the
family to a fiber-preserving family. Since conditions (1) and
(2$^{\prime\prime}$) must remain true in a neighborhood of $t$, we can
cover $D^k$ by convex $k$-cells $B_j$, having corresponding levels
$u_j$ for which (1) and (2$^{\prime\prime}$) hold throughout $B_j$. It
is convenient to rename the $B_j$ so that $u_1<u_2<\ldots<u_r$. Choose
a PL triangulation $\Delta$ of $D^k$ sufficiently fine so that each
$i$-cell lies in at least one of the $B_j$. We will proceed by
increasing induction to deform the family of Klein bottles to be
vertical over the $i$-skeleta of $\Delta$. It will never be necessary
to change the imbeddings over the boundary of $D^k$.

Suppose first that $\tau$ is a 0-simplex of $\Delta$. Let
$j_1<j_2<\ldots <j_s$ be the values of $j$ for which $\tau\subseteq
B_j$.  By $(2^{\prime\prime})$ each intersection circle of $K_\tau$
with each $T_{j_q}$ is isotopic in $T_{j_q}$ to a
fiber of the Seifert fibering. Also, each is an orientation-preserving
simple closed curve in $K_\tau$, so must be isotopic in $M$ to the $a$
loop or the $b^2$ loop in $K_0$. When $m\eq 1$, $b^2$ is the generic
fiber of $M$, and $a$ is not isotopic in $M$ to $b^2$ since $a\eq
b^{2n}$ and $n\neq 1$. When $n\eq 1$, $a$ is the fiber of $M$, and
$b^2$ is not isotopic to $a$ since $a^m\eq b^2$ and $m\neq1$. In
either case, the intersection circles are isotopic in $K_\tau$ to a
fiber on $K_\tau$ (the images of the fibers of $K_0$ under the
imbedding corresponding to~$\tau$). Therefore we may deform the
parameterized family near $\tau$ so that each $K_\tau\cap T_{j_q}$
consists of fibers on $T_{j_q}$ that are images of fibers
of~$K_0$.  Then, using lemma~\ref{fiberpreserving2} successively on
the solid torus $R_{u_s}$, the product regions
$\overline{R_{u_j}-R_{u_{j-1}}}$ for $j\eq j_s,i_{s-1},\ldots,j_2$,
and the twisted $I$-bundle $P_{u_{j_1}}$, deform $K_\tau$ to be
fiber-preserving.

Inductively, suppose that $K_t$ is vertical for each $t$ lying in any
$i$-simplex of $\Delta$. Let $\tau$ be an $(i+1)$-simplex of $\Delta$,
and let $j_1<j_2<\ldots<j_s$ be the values of $j$ for which $\tau$
lies in $B_j$. For each $t\in\partial\tau,\;K_t$ is vertical. By
lemma~\ref{fiberpreserving1} applied to each parameterized family
$K_t\cap T_{j_q}$, we may assume that $K_t\cap T_{j_q}$ consists of
fibers. Again using lemma~\ref{fiberpreserving2} and proceeding from
$R_{u_{j_s}}$ to $P_{u_{j_1}}$, deform the family on $\tau$, keeping
it fixed over $\partial\tau$, to be vertical for all points in
$\tau$. This completes the induction step and the proof of
theorem~\ref{main}.


\begin{thebibliography}{99}
{\footnotesize

\bibitem{ABBJRS} C. Aneziris, A. P. Balachandran, M. Bourdeau, S. Jo,
      T. R. Ramadas, R. Sorkin,  Aspects of spin and statistics
      in generally covariant theories, {\em Int. Jour. Mod. Phys.}
      A4 (1989), 5459-5510.

\bibitem{A} K. Asano, Homeomorphisms of prism manifolds,
    {\em Yokohama Math. J.} 26 (1978), 19-25.

\bibitem{B} F. Bonahon, Diff\'eotopies des espaces lenticulaires,
      {\em Topology} 22 (1983), 305-314.

\bibitem{C-S} S. Cappell and J. Shaneson, Counterexample to the oozing
   problem for closed manifolds, Springer-Verlag Lecture Notes in
   Mathemtics No. 763 (1979), 627-634.

\bibitem{F-S} J. Friedman and R. Sorkin, Spin 1/2 from gravity,
             {\em Phys. Rev. Lett.} 44 (1980), 1100-1103.

\bibitem{F-W} J.  Friedman and  D. Witt, Homotopy is not isotopy for
        homeomorphisms of 3-man\-i\-folds, {\em Topology} 25 (1986), 35-44.

\bibitem{G} A. Gramain, Le type d'homotopie du groupe des diff\'eomorphismes
   d'une surface compacte, {\em Ann. \'Ecole Norm. Sup.} 6 (1973), 53-66.

\bibitem{Giulini} D. Giulini, On the configuration space topology in general
        relativity, Helv. Phys. Acta (1995), 86-111.

\bibitem{G-L} D. Giulini and J. Louko, No-boundary $\theta$-sectors
      in spatially flat quantum cosmology. {\em Phys. Rev. D} 46
      (1992), 4355-4364.

\bibitem{H2} A. Hatcher, Homeomorphisms of sufficiently large
                $\P^2$-irreducible 3-manifolds, {\em Top\-ology}
                15 (1976), 343-347.

\bibitem{H3} A. Hatcher, On  the diffeomorphism group of $S^1\times S^2$,
  {\em Proc. Amer. Math. Soc.} 83 (1981), 427-430.

\bibitem{H1} A. Hatcher, A proof of the Smale Conjecture, $\Diff(S^3)
    \cong \hbox{O}(4),$ {\em Annals of Math.} 117 (1983), 553-607.

\bibitem{I} C. J. Isham, Topological $\theta$-sectors in canonically
      quantized gravity, {\em Phys. Lett. B} 106 (1981), 188-192.

\bibitem{I1} N. Ivanov, Homotopy of spaces of automorphisms of some
            three-dimensional manifolds,
            {\em Soviet Math. Dokl.} 20 (1979), 47-50.

\bibitem{I2} N. Ivanov, Homotopy of spaces of diffeomorphisms of some
    three-dimensional manifolds, {\em J. Soviet Math.} 26 (1984),
     1646-1664.

\bibitem{Jaco} W. Jaco, {\it Lectures on Three-manifold Topology,}
     CBMS Regional Conference Series No.~43 (1977).

\bibitem{Kalliongis-McCullough} J. Kalliongis and D. McCullough,
     Fiber-preserving imbeddings and diffeomorphisms, preprint.

\bibitem{L} F. Laudenbach, Topologie de la dimension trois.
         Homotopie et isotopie, {\em Ast\'erisque} 12 (1974), 1-152.

\bibitem{Orlik} P. Orlik, {\em Seifert Manifolds,} Springer-Verlag
         Lecture Notes in Mathematics Vol. 291, 1972.

\bibitem{Palais} R. Palais, Local triviality of the restriction map for
       imbeddings, {\em Comment. Math. Helv.} 34 (1960), 305-312.

\bibitem{P} R. Palais, Homotopy theory of infinite-dimensional
      manifolds, {\em Topology} 5 (1966), 1-16.

\bibitem{P-R} J. Pitts and J. H. Rubinstein, Applications of minimax
   to minimal surfaces and the topology of 3-manifolds, {\em Proc. Centre
    Math. Analysis} 12, ed. J. Hutchinson and L. Simon, Canberra (1987),
    137-170.

\bibitem{R1} J. H. Rubinstein, One-sided Heegaard splittings of
            3-manifolds, {\em Pacific J. Math.} 76 (1978), 185-200.

\bibitem{R2} J. H. Rubinstein, On $3$-manifolds that have finite
   fundamental group and contain Klein bottles, {\em Trans. Amer. Math.
   Soc,} 251 (1979), 129-137.

\bibitem{R-B} J. H. Rubinstein and J. Birman, One-sided Heegaard
            splittings and homeotopy groups of some 3-manifolds,
            {\em Proc. London Math. Soc.} (3) 49 (1984), 517-536.

\bibitem{Smale} S. Smale, Diffeomorphisms of the 2-sphere, {\it Proc.
     Amer. Math. Soc.} 10 (1959), 621-626.

\bibitem{S3} R. Sorkin, Classical topology and quantum phases: Quantum
      Geons, in {\em Geometrical and algebraic aspects of nonlinear
      field theory,} ed. S. De Filippo, M. Marinaro, G. Marmo, and G.
      Vilasi, Elsevier Science Publishers B.~V. (North
      Holland), 1989.

\bibitem{Waldhausen} F. Waldhausen, Recent results on sufficiently
    large 3-manifolds, {\em Proc. Sympos. Pure Math.} (R. J. Milgram,
    ed.), Vol. 32, Amer. Math. Soc., Providence, RI, 1978, 21-38.

\bibitem{W2} D. Witt, Symmetry groups of state vectors in canonical
     quantum gravity, {\em Jour. Math. Phys.} 27 (1986), 573-592.

\bibitem{Wolf} J. Wolf, {\em Spaces of Constant Curvature,} Publish Or
    Perish Press, 1974.

}
\end{thebibliography}
\end{document}